\documentclass[hein,final,leqno]{shein}
\usepackage{latexsym}
\usepackage{graphics}
\usepackage{epsfig}
\usepackage{epsf}
\usepackage{eucal}
\usepackage{pstcol}
\usepackage{amsfonts}
\usepackage{amstext}
\usepackage{amssymb}
\usepackage{amsmath}
\usepackage{amsbsy}
\usepackage{amsxtra}
\usepackage{amscd}
\newtheorem{ltheorem}{Theorem}           
\newtheorem{lcorollary}{Corollary}[section]  

\newtheorem{lproposition}{Proposition}

\renewcommand{\abstractname}{Abstract.}
\begin{document}
\title{Dual Scattering Channel Schemes Extending the Johns Algorithm}
\author{Steffen Hein
}                     
%
%
\institute{SPINNER GmbH., M\"unchen, Dept TB01
[Numerical Methods], Aiblinger Str. 30,\linebreak DE-83620 Westerham,
Germany\hfill\email{s.hein@spinner.de}}
%
%
\date{}
%
\titlerunning{Dual Cell-Interface Scattering Channel Schemes}
\authorrunning{Steffen Hein}

\maketitle 
%
\begin{abstract}
\vspace{-.2cm}
Dual scattering channel schemes extend the transmission line matrix 
numerical method (\textsc{Johns'} TLM algorithm) in two directions.
For one point, transmission line links are replaced by abstract
scattering channels in terms of paired distributions (characteristic
im\-pe\-dan\-ces are thus neither needed, nor in general defined, e.g.).
In the second place, non-trivial cell interface scattering is admitted
during the connection cycle.
Both extensions open a wide field of applications beyond the range of 
classical time domain schemes, such as \textsc{Yee}'s FDTD method and TLM.
A DSC heat propagation [diffusion] scheme in non-orthogonal mesh,
wherein heat sources are directly coupled to a lossy Maxwell field,
illustrates the approach.
\vspace{5pt} 
\hfill \textbf{MSC-classes}:\,\textnormal{\,65M06,\,78M20,\,80M20}

\vspace{-.2cm}
\end{abstract}

\markboth{{\normalsize \textsc{Steffen Hein}}}
{{\normalsize \textsc{Dual Scattering Channel Schemes}}}

\normalsize
\vspace{-0.3cm}
\begin{quote}\notag
\small
{
\textnormal{Only a man who knows nothing of reason}
\newline
\textnormal{talks of reasoning without strong first principles.}
\begin{flushright}
\textit{Gilbert Keith Chesterton}
\end{flushright}
}
\end{quote}
\vspace{-1.0cm}
\section{Introduction}
It seems a paradox - and is just a typical process in mathematical analysis
that a structure turns simple in a more general setting which at the same time
widens its range of application. Accordingly - ~not very surprising~ - some
ill-famed 'intricacies of the propagator approach to TLM' (sic. Rebel
\cite{diss}, p. 5) virtually vanish if some of its elements are taken as the
building blocks of a more general scheme.
In fact, constructing the latter on essentially these elements in a quasi
axiomatic manner will prove such intricacies to be mere artefacts of an
inadequate framework.\\
The choice of elements proposed in this paper 'generalizes' the Johns
algorithm in two directions.
In the first place, abstract scattering channels replace transmission lines,
which have some unpleasent properties (section \nolinebreak \ref{S:sec2}).
Secondly, non-trivial cell boundary (\emph{interface\/}) scattering
is permitted during the connection cycle.
The schemes thus obtained are characterized by a non trivial two-step
(\emph{connection-reflection\/}) cycle of iteration which exhibits
certain duality relations - ~whence their name.

When P.B. Johns and co-workers introduced the transmission line matrix (TLM)
numerical method in the 1970s \cite{JoBe} it was almost instantaneously
assimilated by the microwave engineering community. In the same audience
the method remained until today subject of assiduous study and extensive
publication. Three conferences explicitely focussed on TLM \cite {A,B,C},
and the monographs of Christopoulos \cite{Ch1} and de Cogan \cite{Co1} deal
in detail with the original ideas as well as with classical applications.

Familiarity with the transmission line picture, and the well known
scattering concept, certainly fostered the acceptance of the TLM method
among microwave engineers.
On the other hand, just so the primary interest turned of course on
applications in their own discipline, rather than onto the
inner algorithmic structure as an object of mathematical analysis.
Over the years still a node more was routinely invented,
with new dispersion characteristics and/or equipped with still another
ingenious stub, designed to model special propagation or transport
phenomena, in varied geometries or boundary conditions.
\cite{Trenkic1} stands somewhat exemplary for this line of research.

Mathematical questions addressing the inner structure of the TLM
algorithm and its potential generalizations have thus apparently
been for a long time of secondary interest.
They have yet not been left completely out of view.
Chen, Ney and Hoefer \cite{Chen} proved equivalence of the original
(expanded) TLM node without stubs \cite{JoBe,Johns0} to the Yee
finite-difference grid \cite{Yee1,KuLue}.
Recently, the non-trivial question of consistence of Johns'
symmetrical condensed node (SCN), cf. Johns \cite{Johns2},
with Maxwell's equations and the intimately related problem
of convergence to a smooth solution for decreasing time step and grid
spacing have been tackled, and in parts solved, by Rebel \cite{diss}.
His thesis presents, by the way, a thorough survey over the ramifications
of TLM until that time (year 2000), without perhaps spending
sufficient attention to its non-orthogonal mesh extensions.

From a quite general viewpoint, viz. widely independent of any
particular physical interpretation, the structure of the stub loaded
(deflected) non-orthogonal TLM algorithm has been analysed in \cite{He5}.
The present paper goes even further and challenges the transmission line
picture at all. The latter, in universally imposing free wave propagation
between cells (with great benefit, at times), induces modeling
limitations under circumstances that are outlined in section
\nolinebreak \ref{S:sec2}.
Many restrictions can be by-passed by replacing transmission lines with
abstract scattering channels in terms of 'paired' distributions.

Dual scattering channel schemes are characterized by a two-step
updating cycle with certain duality relations between the two steps.
The TLM method with its familiar connection-reflection cycle is
trivial \emph{as a dual scheme} in that the connection map
reduces essentially to identity (viz. pure transmission or total
reflection) - \nolinebreak again with modeling limitations.
These can be raised, anew, in permitting non trivial cell interface
scattering during the connection step of iteration.

One major merit of the transmission line method is unconditional
stability ~\cite{Johns3}.
Since this property is usually drawn upon the passivity of linear 
transmission line networks, the question of stability needs a proper
investigation for DSC schemes that do not use lines.
That problem is studied in~\cite{He7}, where it is shown that a wide
class of DSC schemes are in fact unconditional stable.
Due to the convolution type updating scheme
(\emph{Johns' cycle}; cf. section ~\ref{S:sec3}) it is sufficient for
stability that the reflection and connection maps share simple
contraction properties (paraphrased as $\,\alpha$-\emph{passivity} 
in ~\cite{He7}).
In summary, DSC schemes are unconditinally stable under quite general
circumstances, and they are conceptually simple, though a set of
technical definitions is of course ineluctable in a neat theory.
Last but not least, nothing obscure nor 'intricate' should be associated
anymore with the propagator approach.
\vspace{-0.3cm}
\section{Scattering channels}\label{S:sec2}
Any extension of the TLM method that includes heat transfer, fluid flow,
or particle current, for instance, involves scattering channels other than
transmission lines. The latter, for a non vanishing real part of the
characteristic line impedance, inherently impose wave propagation between
cells. Degenerate lines, with a purely imaginary impedance,
still work in diffusion models, cf. \cite{Ch1}, chap.7. Other types of
transport or modes of propagation, such as for example the relativistic
charged particle current treated in \cite{He5}, are very unnaturally and
more or less imperfectly modeled using transmission lines. There is good
reason to get rid of lines in such and other cases within an extended
framework.

A first step towards the definition of more general scattering channels in
TLM has been undertaken in replacing transmission lines with abstract
projections into in- and outgoing field components, cf. \cite{He5}. It was
postulated in this paper that the propagating fields ('link quantities')
allow of a decomposition into a direct sum
\begin{equation}\label{2.1}
z = z_{in} \oplus z_{out} \quad ,
\end{equation}
$z_{in}$ and $z_{out}$ representing the incident and outgoing fields,
respectively. Moreover, it is essential in our understanding of TLM
that the latter have a merely operational meaning in that only the
total field $z$ enters the dynamical model equations (cf. sections
\ref{S:sec3},\ref{S:sec4}).
In singular cases, a physical interpretation can yet still be given
to $z_{in}$, $z_{out}$ on the basis of a special analysis,
e.g. \cite{He4}, Corollary \nolinebreak 2. \\
For the Maxwell field model the technical passage from the transmission
line formulation to the projection operator setting is outlined in
\cite{He5}, Appendix ~A.

In the classical TLM setup the connection map simply transfers
without further modifications the quantities outgoing from a cell
into quantities incident at neigbouring cells or rejected at some
totaly reflecting electric or magnetic wall.
This is in perfect harmony with the behaviour of a propagating
electromagnetic field, the components of which are tangential to the cell
boundary (as the link quantities always are in a classical TLM cell,
cf. \cite{He1}) and that is thus not subject to refractive scattering
at the cell face, even if the medium changes there.

The situation is clearly not thus simple for arbitrary propagating quantities.
To circumvent any modeling restrictions, the connection cycle of a
non-trivial DSC scheme comprises cell interface scattering from the outset.
Nodal and cell face scattering thus enter a kind of duality relation
that becomes visible, for instance, in an apparent symmetry of the model
equations in their most general form (\ref{3.16}, \ref{4.3n}, \ref{4.3p}).
Nodal and cell boundary scattering may in fact be of equal importance and
sometimes boundary scattering plays even the leading r\^ole in a DSC
algorithm.

In the generalized setup, just as in the traditional TLM framework,
scattering channels interconnect a {\it node\/}, viz. a suitably defined
{\it centre\/} of a mesh cell, with {\it ports\/} at the cell boundary.
The channels are yet no longer represented by transmission lines.
With respect to a computed physical field in D-dimensional configuration
space, they simply form a pair of scalar or vector valued distributions,
transposed over a distance in space, which test the field within the
cell and on its boundary.
A DSC scattering channel will thus be defined, precisely, as a pair of
continuous linear functions (${\,p \,, \; p \sptilde\,}$) which act
on a class of (suitably smooth real or complex) vector fields $Z$
in configuration space, such that $p$ has its support on a cell
face and ${p \sptilde}$ is connected to $p$ via pull back into the node,
i.e.: Given any notions of \emph{centre} of cell and face, as well as
the spatial translation
${s: \mathbb{R}^D \rightarrow \mathbb{R}^D}$
that shifts the centre of a cell (\emph{node\/}) into the centre
of the face where $p$ has its support, then the {\it nodal image}
$p\sptilde$ of $p$ is defined as the distribution
\begin{equation}\label{2.2}
(\, p\sptilde \, , \; Z \, ) \, :\, = \, ( \, p \circ s \, , \; Z \,) \,
= \, ( \, p \,, \; Z \circ s^{-1} \, ) \; ,
\end{equation}
and the pair (${ \, p,\; p\sptilde\,}$) is called a {\it scattering channel}.
Equivalently, a scattering channel can of course be identified with
(${\, p \, , \; s\,}$) or even simply with the \emph{port p\/}, the pertinent
shift and nodal image then being tacitly understood. 
The concept should in fact not be handled in too rigid a fashion - ~ and
there is no need to do so.
In certain applications the support of the port distribution may be
extended over a neighbourhood of a face, or the node distribution may
rather be thought of as a mean over the entire cell (in the way familiar
from finite volume methods).
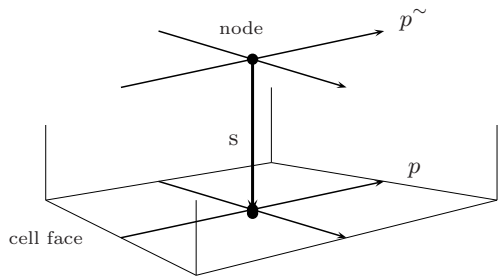
\begin{figure}[!h]\centering
\vspace{-1.0cm}
\setlength{\unitlength}{1.cm}
\begin{pspicture}(0.0,0.0)(6.0,4.5)\centering
\psset{xunit=1.0cm,yunit=1.0cm}
\psline[linewidth=0.1mm]{-}(0.0,1.0)(2.0,0.0)
\psline[linewidth=0.1mm]{-}(2.0,0.0)(6.0,1.0)
\psline[linewidth=0.1mm]{-}(6.0,1.0)(3.0,1.5)
\psline[linewidth=0.1mm]{-}(3.0,1.5)(0.0,1.0)
\psline[linewidth=0.1mm]{-}(0.0,1.0)(0.0,2.0)
\psline[linewidth=0.1mm]{-}(2.0,0.0)(2.0,1.0)
\psline[linewidth=0.1mm]{-}(6.0,1.0)(6.0,2.0)
\psline[linewidth=0.1mm]{-}(3.0,1.5)(3.0,2.5)
\psline[linewidth=0.2mm]{->}(1.0,0.5)(4.5,1.25)
\psline[linewidth=0.2mm]{->}(1.5,1.25)(4.0,0.5)
\psline[linewidth=0.2mm]{->}(1.0,2.5)(4.5,3.25)
\psline[linewidth=0.2mm]{->}(1.5,3.25)(4.0,2.5)
\psline[showpoints=true,linewidth=0.4mm]{->}(2.75,2.875)(2.75,0.825)
\psline[showpoints=true,linewidth=0.4mm]{-}(2.75,0.875)(2.75,0.825)
\rput(2.6,3.3){\scriptsize node}
\rput(4.9,3.4){\small $p \sptilde$}
\rput(2.5,1.8){\small s}
\rput(0.0,0.5){\scriptsize cell face}
\rput(4.9,1.4){\small $p$}
\end{pspicture}
\caption{\textsl{Ports on a cell face with their nodal images.}}\label{F:1}
\end{figure}

Needless to say that the stressed duality between nodal and cell boundary
scattering is not to be misunderstood in the narrow sense of category
theory.
Here, it refers simply to the observation that a set of propositions are
valid, modulo symmetry in certain terms, in the two scattering situations
- \nolinebreak which of course reflects the paired distribution concept
of scattering channel and the already mentioned symmetry of the pertinent
model equations in their most general form.
A parallel symmetry then clearly characterizes the structure of the 
reflection and connection maps that solve these equations.

Cell boundary scattering is, for the rest, not thus new an option:
Already in the TLM model for superconducting boundary \cite {He0}
cell face s-parameters and boundary stubs have been introduced for
solving the discretized London equations, cf. also \cite{He4}.

Despite the abolition of transmission lines, by their replacement
with abstract scattering channels the computed ('physical') fields
can still be represented -- in the way familiar from the classical
TLM method -- as sums of \emph{in- and outgoing} scalar or vector fields
\begin{equation}\label{2.3}\centering
z \; \; = \; \; z_{in} \, + \, z_{out} \quad . 
\end{equation}
No physical interpretation or propagation property is, however, in general
ascribed to $z_{in,\, out}$.
In fact, these quantities are merely operationally defined by means of the well
known \emph{Johns cycle} of iteration
\begin{equation}\label{2.4}\centering
\begin{split}
\begin{aligned}
z_{in} \; &= \; ( \; \mathcal{C} [z_{out}] + e \; )&&\longleftarrow && \\
z^{p} \; &= \; z_{in} + z_{out} && &&\; \, \uparrow \\
&\, \downarrow \; t + \tau/2 && &&t + \tau \\
z_{out} \; &= \; \mathcal{R} [z_{in}] && &&\; \, \uparrow \\
z^{n} \; &= \; z_{in} + z_{out} &&\longrightarrow && \qquad \, .
\end{aligned}
\end{split}
\end{equation}
$\mathcal{R}$ and $\mathcal{C}$ denote the node and cell face propagators
(or so-called reflection and connection maps - the latter including
now cell boundary scattering, and $e = e(t)$ induces any excitation.
Note again that $z_{in}$, $z_{out}\,$ are so far \emph{purely operational}
quantities,
i.e. only the total fields $z$ enter the model equations, while $z_{in}$,
$z_{out}$ are {\it in general\/} bare of any physical meaning (a physical
interpretation in terms of an energy flow still exists within the classical
Maxwell field TLM model, cf. \cite{He4}, Corollary \nolinebreak 2).

As will be seen in the next section, the structures of the propagators
${\,\mathcal{R}\,}$ and ${\,\mathcal{C}\,}$ are very similar in the
general DSC scheme, thus reflecting the dual r\^ole that nodal and
boundary scattering play therein.
In a sense, precised in section \nolinebreak \ref{S:sec3},
${\,\mathcal{R}\,}$ and ${\,\mathcal{C}\,}$ are the discrete
\emph{convolution integrals\/}
that in every Johns cycle strictly solve the model equations.
In fact, the Johns cycle can be looked at as basically a two-step
convolution method for solving certain types of explicit finite
difference equations in time.

Note that the model equations can in principle be directly solved
inasmuch as they provide complete recurrence relations, cf.
sections \nolinebreak \ref{S:sec3}, \ref{S:sec4}.
The scattering approach (using Johns' cycle of convolutions) offers,
however, important advantages. Thus, it provides clear cut reliable
stability criteria that make the DSC algorithm unconditionally stable
under very general circumstances ~\cite{He7}.
\vspace{-0.3cm}
\section{The elements characterizing DSC schemes}\label{S:sec3}
So far, we dealt on a largely informal level with some typical
traits of the TLM algorithm that either characterize DSC schemes
in like manner, or which have to be modified in a specified way
in order to attain a greater generality.
We are still bound to keep our introductory promise and give a
coherent description of DSC schemes in terms of some quasi axioms
that condense their distinctive properties.
Of course, we shall not really pursue axiomatics, here, in the sense
of building a new theory on a complete set of basic assumptions.
(Nor are we adopting a dogmatic attitude and going to fix a rigid
framework that with certainty must be modified sooner or later in
order to face a particular problem.)
The emphasis is rather on compiling on a largely preliminary basis
those formal elements that, in essence, lead to the peculiar structure of
DSC schemes (in general), and of the TLM method (in particular), without
being distracted by unnecessary information, such as mesh topology
and geometry, s-parameters, etc., which characterize only a singular
physical interpretation. A basic set of technical notions need simply
to be clarified. \\
In the following, '{\it simple}' definitions are visualized in writing
the defined object(s) in italics, more crucial or technical ones are
explicitely designated as \textbf{Definition}.

Until further notice, a {\it mesh\/} denotes a non-void finite family
(i.e. an indexed set) of elements named {\it cells\/},
which are sets in their turn, and share the following properties.
Each cell $\zeta$ contains an element $n_\zeta\,$, called its
{\it node\/}, and a finite family
$\partial\zeta = \{ \partial \zeta_{\, \iota} \}\/$,
called the (cell) {\it boundary}.
The latter is built up of elements ${\partial \zeta_{\, \iota}\,}$,
named {\it faces}, which are sometimes simply written $\iota$ in the
place of ${\partial \zeta_{\, \iota}}$.

\begin{definition}\label{D:1}
A mesh $M$ is called \emph{regular\/}, if and only if it satisfies the
following requirements of \emph{simplicity} (S)
and \emph{connectedness} (C):
\begin{itemize}
\item[(S)]
   Every node belongs to exactly one cell and every face to at most two
   cells in $M$.
\item[(C)]
   For every two cells $\zeta^i, \zeta^j \in M$, there exists a
   connecting sequence \newline
   ${s = (\zeta^\kappa )^k_{\kappa = 0} \in
   M^\mathbb{N}}$, such that $\zeta^0 = \zeta^i, \zeta^k =
   \zeta^j$ and every two subsequent cells $\zeta^\kappa,
   \zeta^{\kappa +1}$ in $s$ have a common face, for $0 \leq \kappa < k$.
\end{itemize}
\end{definition}

Also - certainly not too misleading:
any two cells with a common face are called {\it adjacent\/} or
{\it neighbouring\/} cells, and the common face a {\it connecting face\/}
or {\it interface}. \\
By a first postulate,
\emph{DSC meshes are always regular meshes}.

The {\it state space\/} is any product of real or complex normed linear
spaces labelled by mesh cells
\begin{equation}\label{3.1} 
\mathcal{S}= \prod \nolimits_{\, \zeta \in M} \, \mathcal{S}_{\zeta} \quad . 
\end{equation}
In addition, we require that a DSC state space always contains a non-void
subspace $\mathcal{P} \subset \mathcal{S}\,$,
named the space of {\it propagating fields\,},
which on every $\mathcal{S}_\zeta$ reduces to a product of 'squared' 
spaces in the following precise sense
\begin{equation}\label{3.2}
\mathcal{P}_\zeta = \mathcal{P}\cap\mathcal{S}_{\zeta} = \prod
\nolimits_{\,\iota\,\in\,\partial\zeta} \;
P_{\,\iota,\,\zeta}^{\,2}\quad .
\end{equation}
In less formal language: every propagating field $z$ splits over the
cells into a sequence of pairs
\begin{equation}\label{3.3}
z_\zeta = (\,z_\iota,\,z_\iota\sptilde\,)_{\iota \,\in\,\partial\zeta}
\end{equation}
labelled by the faces of the cell boundary. In other words, there is a
canonical automorphism of normed linear spaces on $\mathcal{P}$, which on
every $\mathcal{P}_{\zeta}$ reduces to
\begin{equation}\label{3.4}\centering
\begin{split}
\begin{aligned}
nb \; \; : &\quad \mathcal{P}_{\zeta} &&\longrightarrow
&&\quad \mathcal{P}_{\zeta} \\
&(\,z_\iota\,,\,z_\iota\sptilde\,) &&\longmapsto
&&(\,z_\iota\sptilde\,,\,z_\iota\,) \quad .
\end{aligned}
\end{split}
\end{equation}
$nb$ is obviously involutary (${\;nb^{\,2} = Id\;}$), and is
called the {\it node-boundary map}.

The components $z_\iota$ and $z_\iota \sptilde$
in \eqref{3.3}, \eqref{3.4} are named the
{\it port\/} (or {\it face\/}) {\it component}, and the
{\it node component}, respectively, of
$z \, = \, (\, z_{\iota} \,, \, z_{\iota}\sptilde\,$).
They are usually written $z^p = z_{\iota} = \pi^{p} ( z )$ and
$z^n = z_{\iota} \sptilde = \pi^{n} ( z )\;$ with projections
$\pi^{p}\,$, $\pi^{n}\,$
that are canonically extended over the entire space $\mathcal{P}$.

Let ${J :\,= \bigcup_{\zeta \in M} \partial \zeta}$ be the set of all
faces in $M$ (remember that $\partial \zeta$ is defined as a union
of faces). Then, in virtue of \eqref{3.2}, $\mathcal{P}$
splits completely into subspaces
\begin{equation}\label{3.5}
\quad \quad \mathcal{P} = \prod \nolimits_{\, \iota \, \in \, J}
\mathcal{P}_{\iota} \qquad \text{with} \qquad
\mathcal{P}_{\iota} \, :\,= \,
\prod \nolimits_{\, \zeta \; ; \; \iota \, \in \, \partial \zeta}
\; P_{\, \iota , \, \zeta}^{\,2} \quad .
\end{equation}
A {\it DSC process} is a step function of time
\begin{equation}\label{3.6} 
pr:[\, 0 \,, \, T\, ) \longrightarrow \mathcal{S} \quad ,
\end{equation}
such that ${\pi^p \circ pr(t)}$ and ${\pi^n \circ pr(t-\tau/2)}$ are
constant on every time interval ${[\mu\tau,(\mu+1)\tau)}$, ${\mu \in
\mathbb{N}\,}$, where they are defined. In other words,
port components of a DSC process switch at integer multiples of the
time step $\tau$ while node quantities switch at odd integer multiples
of $\tau /2$.

Given a process, a state $z$ with its entire history up to time $t$ is
usually written as a 'back in time running' sequence
\begin{equation}\label{3.7}
[\,z\,](\,t\,) \;
:\,= \; (\, z \, (\,t- \mu \tau\,)\,)_{\mu \in \mathbb{N}} \quad ,
\end{equation}
expanding so the domain of definition eventually to the negative time axis
in the trivial way, i.e. $z(s)\, :\,= \, 0$ \, for $s<0$.

By this convention, we assign thus to index $\mu$ the (varying) state back
in the past from present time $t$
\begin{equation}\label{3.8} 
[\,z\,]_\mu (\,t\,) \; :\,= \; z \, (\,t- \mu \tau\,) \quad , \quad
\end{equation}
rather than the (fixed) state $z(\mu \tau)$ - which has the technical
advantage that $\mu$ so is directly related to a {\it time difference\/}
(and eventually to an order of a finite difference equation in time),
rather than to an {\it absolute time\/} (which is quite uninteresting,
in general).

Functions defined on back in time running sequences, such as \eqref{3.7},
are called {\it causal} functions (or {\it propagators\/}). Any such map
is a discrete analogue to a causal Green's function integral, as has
already been outlined in \cite{He4}.

Let $X_{0}^{\mathbb{N}}$ denote the set of all sequences with an arbitrary,
but finite, number of non-vanishing elements in a linear space $X$.
For every mesh cell $\zeta$ and face $\iota$ in $M$ consider then
the subspaces of propagating fields
$\mathcal{P}_{\zeta}^n \, :\,= \, \pi^{n} (\mathcal{P}_\zeta)\,$,
$\mathcal{P}_{\iota}^p \, :\,= \, \pi^{p} (\mathcal{P}_\iota)$.

\begin{definition}\label{D:2}
A \emph{reflection map} (in $\zeta \in M\,$) is a
(possibly time dependent) causal operator
\begin{equation}\label{3.9}\centering
\begin{split}
\begin{aligned}
\mathcal{R}_{\zeta} \; :\; (\mathcal{P}_{\zeta}^{n})_{0}^{\mathbb{N}} \;
&\longrightarrow \; \mathcal{P}_{\zeta}^{n} \\
[z^n] \; \; \, &\longmapsto \; \mathcal{R}_{\zeta} [z^n] \;\;,
\end{aligned}
\end{split}
\end{equation}
and a \emph{connection map} (in $\iota \in J\,$) is a
(likewise possibly time dependent) causal operator
\begin{equation}\label{3.10}\centering
\begin{split}
\begin{aligned}
\mathcal{C}_{\iota} \; : \; (\mathcal{P}_{\iota}^{p})_{0}^{\mathbb{N}} \;
&\longrightarrow \; \mathcal{P}_{\iota}^{p} \\
[z^p] \; \; \, &\longmapsto \; \mathcal{C}_{\iota} [z^p] \;\;.
\end{aligned}
\end{split}
\end{equation}
\end{definition}

\hspace{-15pt}
Also, a {\it DSC system over M} is a pair ($\mathcal {C}$,
$\mathcal {R}$)
consisting of any two families
\begin{equation}\label{3.11}\centering
\; \; \mathcal{C} = \{\mathcal{C}_{\iota}\}_{\iota \in J}
\quad \text{and} \quad
\mathcal{R} = \{\mathcal{R}_{\zeta}\}_{\zeta \in M} \quad
\end{equation}
of connection and reflection maps.
Note that sometimes also the \emph{families} are called the connection
and reflection map. For the sake of algorithm stability,
these maps (in whatever meaning) should share certain contraction properties
studied in ~\cite{He7}. 

An {\it excitation\/} is only a distinguished process with values in
the mesh boundary states. More precisely, let
${B :\,= \{ \iota \; | \; \text{$\iota \in J$ and $\iota$
is not an interface} \} }$
be defined as the {\it mesh boundary}, then an {\it excitation\/} 
is a process
\begin{equation}\label{3.12}\centering
\begin{split}
\begin{aligned}
e \; : \; [\,0, \, T) \; &\longrightarrow \;
{\prod \nolimits}_{\, \iota \, \in \, B}
\mathcal{P}_{\iota} \\
\qquad t \; &\longmapsto \; e(t)=\pi^p\circ e(t) \; ,
\end{aligned}
\end{split}
\end{equation}
i.e. $e$ is a port process, and hence switches at entire multiples
of the time step $\tau$, and $e$ generates \emph{non-interface\/}
(\emph{mesh boundary\/}) states.

\begin{definition}\label{D:3}
The \emph{DSC process generated by} ($\mathcal{C}$, $\mathcal{R}$)
\emph{and excited by} $e\,$ is the unique process
$z(t)=(z^p, z^n)(t)\/$ which at every time $t \in [0,T)$ satisfies
\begin{equation}\label{3.13}\centering
\begin{split}
\begin{aligned}
z^{\,p}(t)\; &\,=\; \; nb\circ z_{in}^{\,n}(t+\tau/2)+z_{out}^{\,p}(t) \\
z^{\,n}(t)\; &\,=\; \; z_{in}^{\,n}(t)+nb\circ z_{out}^{\,p}(t+\tau /2)
\quad ,
\end{aligned}
\end{split}
\end{equation}
the right-hand side being recursively defined through cyclic iteration of 
\begin{equation}\label{3.14}\centering
\begin{split}
\begin{aligned}
z_{in}^{\,n}(t+\tau/2) \; \; &:\,= \; \; nb\circ
[ \, \mathcal{C}[z_{out}^{\,p}](t)+e(t) \, ] \\
z_{out}^{\,p}(t+\tau) \; \; &:\,= \; \; nb\circ
\mathcal{R}[z_{in}^{\,n}](t+\tau /2) \\
t\; \; &:= \; \; \, t \, + \, \tau
\end{aligned}
\end{split}
\end{equation}
(in that order) with, initially, $z_{in}^{\,n}(0)=z_{out}^{\,p}(0)=0$. \\
Remember that $nb$ denotes the node-boundary map \eqref{3.4}.
\end{definition}

In \eqref{3.14} $\mathcal{C}$ and $\mathcal{R}$ stand of course for
application of all propagators
$\mathcal{C}_{\iota}$ and $\mathcal{R}_{\zeta}$
in the pertinent families (over $J$ and $M$, respectively). Note that
the order of application within the families is unimportant in virtue
of the pairwise disjointness of all $\mathcal{P}_{\iota}$, and of all
$\mathcal{P}_{\zeta}$ - ~which obviously implies that either
$\mathcal{C}$ and $\mathcal{R}$ are \emph{completely parallelizable}
at every time step.

It follows immediately that $z^n$ and $z_{in}^{n}$
thus defined are node processes, hence switch at odd integer multiples
of $\tau/2$, while $z^p$, $z_{out}^{p}$ are port processes
(so they carry their superscripts aright).\\
Equations \eqref{3.13}, \eqref{3.14} can still be simplified to
\begin{equation}\label{3.13'}\centering
\begin{split}
\begin{aligned}
\qquad z^{\,p,\,n} \; \; &= \; \; z_{in}^{\,p,\,n}+z_{out}^{\,p,\,n}
\qquad \qquad \qquad \text{and} \\ 
\qquad \quad z_{in}^{\,p}(t) \; \; 
&= \; \; \mathcal{C}[z_{out}^{\,p}](t) + e(t) \\
\qquad \qquad z_{out}^{\,n}(t+\tau /2) \; \; 
&= \; \; \mathcal{R}[z_{in}^{\,n}](t+\tau /2) \quad ,
\end{aligned}
\end{split}
\end{equation}
writing as usual $z_{in}^{\,p}$ and $z_{out}^{\,n}\/$ for
\begin{equation}\label{3.15}\centering
\begin{split}
\begin{aligned}
z_{in}^{\,p}(t) \; \; &:\,= \; \; nb\circ z_{in}^{\,n}(t+\tau /2) \\
nb\circ z_{out}^{\,p}(t) \; \; 
&=\,: \; \; z_{out}^{\,n}(t-\tau /2) \quad .
\end{aligned}
\end{split}
\end{equation}
Comparing this to the TLM usage we note that in \cite{He5} port
quantities $z_{in}^{p}$, $z_{out}^{p}$ are first introduced. With
these are then node quantities $z_{in}^{n}$, $z_{out}^{n}$ identified
(modulo the time shifts $\pm \tau/2$, just as in \eqref{3.14}) without
yet explicitely mentioning the node-boundary isomorphism $nb$. We shall 
sometimes follow this usage and omit the symbol $nb$ where this cannot
lead to confusion.

Nothing has been said, so far, about physical interpretations or any
implemented dynamical equations. In fact, the characteristic structure of
the DSC algorithm is entirely laid down with the given definitions.
As will be seen in the next section, typical features and facts, some quite
familiar from TLM, are derived straight away with only the above elements.

With respect to the dynamical \emph{model equations} - \nolinebreak which
the algorithm has ultimately to solve and that determine the propagators
${\mathcal{R}_{\zeta}}$, ${\mathcal{C}_{\iota}}$,
cf. section \nolinebreak \ref{S:sec4} \nolinebreak -
we reiterate the important general agreement that only {\it total fields\/}
$z^p\/$, $z^n$, {\it not\/}, however, their {\it incident\/} and
{\it outgoing components separately\/}, shall enter these equations.
Accordingly, we consider only model equations between total fields.
Quite generally, and modulo further restrictions (inferred in the next
section), the DSC model equations should be of the types
\begin{equation}\label{3.16}\centering
\begin{split}
\begin{aligned}
\mathcal{F}^{n}[z_{+}^{n}][z^{p}] \; \; \; &\equiv &&\, \, 0 \\
\mathcal{F}^{p}[z_{+}^{p}][z^{n}] \; \; \; &\equiv &&\, \, 0 \quad ,
\end{aligned}
\end{split}
\end{equation}
with causal functions $\mathcal{F}^{n}$, $\mathcal{F}^{p}$ and shortly
${[\,z_{\pm}\,]}$ \,${:\,=}$ \, ${[\,z\, ](\, t \pm \tau /2\,)}$.
The ${\tau/2}$ time shifts synchronize node and cell boundary switching
in \eqref{3.16}, such that the equations can be {\it strictly\/} solved,
for {\it every} time ${t \, \in [\, 0\, ,\, T\, )}$ and are well-posed,
in this sense.
Of course, time shifts by ${-\tau/2}$ would also lead to synchronization.
The resulting equations would, however, conflict with the causality
property of $\mathcal{R}$ and $\mathcal{C}$.

Inspection of equations \eqref{3.16} shows that $\mathcal{F}^{n}$ affects
only the reflection cycle, while $\mathcal{F}^{p}$ has impact only on the
connection cycle. We are now dealing with the model equations in some more
detail.
\vspace{-0.3cm}
\section{The model equations}\label{S:sec4}
The physical interpretation of a DSC system fixes, intuitively speaking,
the terms in that states in $\mathcal{P}$ are read as physical fields.
More deliberately, certain states in $\mathcal{P}$ are interpreted as
distributional values (finite integrals, e.g.) of physical fields, which
are localized in a mesh cell system.

Any interpretation requires, hence, in the first instance a geometric
realization of the underlying regular mesh, wherein the relations
between abstract cells, nodes, and boundary faces which characterize $M$
are translated into relations between geometric objects, bounded subsets
of ${\mathbb{R}^{D}}$, such as (in general) polyhedral mesh cells with
their faces, e.g.

Given any geometric realization of $M$, a \emph{physical interpretation\/}
of a DSC system over $M$ is, precisely, a family 
${I = \{ I_{\zeta}^{\iota} \}_{\zeta \in M,\; \iota \in \partial \zeta }}$
of continuous linear functions
\begin{equation}\label{4.1}
I_{\zeta}^{\iota} : \mathcal{E}_{\zeta}^{\iota} 
\longrightarrow \mathcal{P}_{\zeta} \quad ,
\end{equation}
each defined on a space ${\mathcal{E}_{\zeta}^{\iota}}$ of smooth
m-component
vector fields in D-dimensional configuration space
(\,i.e. ${\mathcal{E}_{\zeta}^{\iota}
\subset C^{\infty}(\mathbb{R}^{D})^{m}}$;
${m\in \mathbb{N}}$ depending on ${(\iota,\zeta)\;}$),
such that ${I_{\zeta}^{\iota}}$
has its distributional support on a cell face and its range in
${\mathcal{P}_{\zeta}^{p} = \pi^{p}(\mathcal{P}_\zeta)}$.

Note that index $\iota$ in \eqref{4.1} may optionally be read as a cell
face label {\it or} as a multiindex referring to a set of ports on
the same face. 
Since we are dealing with vector-valued distributions
(with range in $\mathcal{P}\,$) and the
support of every ${I_{\zeta}^{\iota}}$ is required to be \emph{localized
on} a cell face (which can be weakened to at least \emph{associated to} 
a face), there is essentially no difference in reading ${z_{\iota,\zeta}^{p}
= I_{\zeta}^{\iota}(Z)}$ as a cell face state vector, or as an array of
components (labelled by port indices) of such a vector.
Thus, index $\iota$ in \eqref{4.1} may be thought of as implicitely
labelling a subset of ports on face ${\iota \in \partial \zeta}$.

Attention is also drawn to the fact that the
functions $I_{\zeta}^{\iota}$ are \emph{not} required to be
\emph{surjective} onto ${\mathcal{P}_{\zeta}^{p}\,}$ (i.e. not every state
in ${\, \mathcal{P}_{\zeta}^{p}\,}$ must be directly related do a
distribution in $I\,$).
There is, for instance, no need to exclude from
${\,\mathcal{P}_{\zeta}^{p}\,}$
any function or linear combination of fields in \emph{different} spaces
${\,\mathcal{E}^{\iota}\,}$ (which may represent a spatial finite difference
of fields, as in the approximate gradient of our sample model in section~
\ref{S:sec5}\,).

The evaluation of nodal fields goes quasi pick-a-pack with ${\it I}$ by
applying the scattering channel concept of section~\ref{S:sec3}.

Let ${s_{\zeta}^{\iota}}$ denote the translational shift in
${\mathbb{R}^{D}}$ from any node ${n_{\zeta}}$ into the (centre of the)
face where ${I_{\zeta}^{\iota}}$ has its support, cf. fig \nolinebreak
\ref{F:1}, and assume without loss of generality that $S_{\zeta}^{\iota} :
Z(x) \mapsto Z(x+s_{\zeta}^{\iota})$ is an inner map in
$\mathcal{E}_{\zeta}^{\iota}$
(otherwise take the closure of $\mathcal{E}_{\zeta}^{\iota}$ under such
transformations).
Then clearly holds
\begin{lproposition}{\label{P:1}}
For every $I_{\zeta}^{\iota} \in I$ there exists exactly one function
${I_{\zeta}^{\iota \sptilde}:\mathcal{E}_{\zeta}^{\iota}\rightarrow} \\
{\mathcal{P}_{\zeta}^{n}} \,\, ( \, $note \, $I_{\zeta}^{\iota \sptilde}
\notin I \, $ ),
such that the following diagram is commutative

\begin{equation}\label{4.p1}\centering
\begin{CD}
\mathcal{P}_{\zeta}^{p} @> nb >> \mathcal{P}_{\zeta}^{n} \\
@AI_{\zeta}^{\iota}\quad AA @AA \quad I_{\zeta}^{\iota \sptilde} A \\
\mathcal{E}_{\zeta}^{\iota} @>> S_{\zeta}^{\iota} > \mathcal{E}_{\zeta}^{\iota}
\end{CD}
\end{equation}
i.e. ${I_{\zeta}^{\iota \sptilde}
\circ S_{\zeta}^{\iota} = nb 
\circ I_{\zeta}^{\iota}}$.
\end{lproposition}
\begin{proof}
Mere retrospection of definitions.
\end{proof}

In the terminology of section~\ref{S:sec2} is
${\; I_{\zeta}^{\iota \sptilde}}$
the nodal image of the port(s)
${I_{\zeta}^{\iota}}$,
and the pair of distributions
( ${I_{\zeta}^{\iota} \; , \; nb \circ I_{\zeta}^{\iota \sptilde}}$ )
forms a scattering channel.

The dynamical DSC \emph{model\/} 
\emph{equations\/} are,
in the line of the preceeding, to be read as
finite difference equations \emph{in time} between states \newline
${z=(z^p_{\iota , \zeta}, z^n_{\iota , \zeta} ) \in \mathcal{P} \;}$
that have an interpretation as distributional values of physical fields
${Z \in \mathcal{E}_{\zeta}^{\iota}\;}$,
\begin{equation}\label{4.2}
z^p_{\iota, \zeta}=I^\iota_\zeta (Z)\quad , 
\quad z^n_{\iota, \zeta}=I^{\iota \sptilde}_\zeta (Z) \quad .
\end{equation}
Unlike classical FD equations between pointwise evaluated physical fields,
the DSC model equations interrelate in many cases finite integrals over
a line segment or face, e.g. - pointwise evaluation (with a Dirac measure
as distribution) not excluded.
The DSC approach is, in this respect, by far more versatile than the
classical finite difference time domain method.
So, the former permits, for instance, of generalizing to a non-orthogonal
mesh Johns' TLM method \cite{Johns2} in much a simpler way, cf. \cite{He2},
than the FDTD approach allows for Yee's method of approximation to
Maxwell's equations \nolinebreak \cite{Yee1}.

A central principle underlying DSC schemes is
{\it near-field interaction}.
Like causality, this is already implicit in the (domains of) definition
of the reflection and connection maps.
Near-field interaction simply spells that only the fields in the
immediate neighbourhood of a state $z\,$ - \nolinebreak precisely
only those in ${\mathcal{P}_{\zeta}\/}$ if
${z \, \in \, \mathcal{P}_{\zeta}^{n}}$, and those in
${\mathcal{P}_{\iota}\/}$ if ${z \, \in  \, \mathcal{P}_{\iota}^{p}}$,
cf. (\ref{3.2}, \ref{3.5}), along with their history, of course
\nolinebreak -
determine the evolution of that state on the next updating step
(note, this refers to fields evaluated in $\mathcal{P}$ and not,
for instance, to an exterior potential, which may still induce
a time dependence of $\mathcal{R}$ or $\mathcal{C}\,$).

In other words, an updated nodal state depends only on states of the
pertinent cell, including its boundary, while the evolution of a port
state is determined by states of the respective face and by nodal
states of the adjacent cells.

It follows that the model equations \eqref{3.16}
split into the two families
\begin{equation}
\begin{aligned} \label{4.3n}
&\mathcal{F}_{\zeta}^{n}[z_{\zeta +}^{n}][z_{\zeta}^{p}] 
\, \, &&\equiv \, \, &&0 \; , \; &
&z_{\zeta} \in \mathcal{P}_{\zeta} \; , \; &&\zeta \in M
\end{aligned}
\end{equation}
\begin{equation}
\begin{aligned} \label{4.3p}
&\mathcal{F}_{\iota}^{p}[z_{\iota +}^{p}][z_{\iota}^{n}]
\, \, &&\equiv \, \, &&0 \; , \; &
&z_{\iota} \in \mathcal{P}_{\iota} \; , \; &&\iota \in J \; .
\end{aligned}
\end{equation}
(Remaining aware of the dependence upon $\zeta$ and $\iota$ of these
equations, we can in general omit the subscripts, if there is no danger
of confusion.)

Since the following analysis runs perfectly parallel for the
dual equations, from now on it is confined to the implications of
\eqref{4.3n}.
- \nolinebreak The reader may write down
parallel statements for dual equations, at times, by exchanging
port for nodal and incident for outgoing quantities,
starting, for instance, with the cell boundary version
(with $\mathcal{C}$ in the place of $\mathcal{R}$) of the
following:

\begin{definition}\label{D:4}
Let ${I \, = \, [\, 0, \, T )\/}$ be a \emph{finite} intervall
(\,i.e. $T \in \mathbb{R}_{+} \,$).
Then we shall say that $\mathcal{R}$ \emph{generates
solutions} of the model equations \eqref{4.3n} \emph{on $I$},
if and only if for every sequence of incident nodal fields
$[\,z_{in}^{n}\,]$ the (obviously unique) process \;
${z \, = \, z_{in} \, + \, z_{out}}$\; which is recursively given by
\begin{equation}
\begin{split}
\begin{aligned}\label{4.d4}
z_{in}^{p}( t \, - \, \tau/2 ) \quad &= \quad nb \, z_{in}^{n}( \, t \, ) \\
z_{out}^{n}( \, t \, ) \quad 
&= \quad \mathcal{R}\, [ \, z_{in}^{n} \, ]\, (\, t \, ) \\
z_{out}^{p}( t + \tau/2 ) \quad 
&= \quad nb \, z_{out}^{n}( \, t \, ) \quad
\end{aligned}
\end{split}
\end{equation}
algebraically solves equations \eqref{4.3n} (identically on $I$).
Sometimes, we then simply say that
$\mathcal{R}$ \emph{solves} the model equations \emph{on that\/} 
(\emph{finite\/} !) \emph{interval}.
\end{definition}

\begin{remark}
\, \hfill
\begin{itemize}
\item[(i)]
Note that Definition \ref{D:4} refers to a purely
\emph{algebraic property} of $\mathcal{R}$ that is not yet related
to any stability questions, for instance
(cf. also the Remark to Theorem \nolinebreak \ref{T:1}\,).
\item[(ii)]
If $\mathcal{R}$ solves equations \eqref{4.3n} on $I$, then in particular
every process generated by ($\mathcal{C}$, $\mathcal{R}\/$) in the sense of
Definition \ref{D:4} solves \eqref{4.3n} on $I$, since every such process
obviously satisfies \eqref{4.d4} (cf. \eqref{3.15}).
\end{itemize}
\end{remark}

The least awkward (yet fortunately frequently encountered) situation
is \\ brought abount with homogeneous linear equations, i.e. for the evolution
of nodal states 
\begin{equation}\label{4.4}\centering
\mathcal{F}^{n} [z_{+}^{n}] [z^{p}] = \sum \nolimits_{\mu=0}^{\infty} \phi_{\mu}
z^{n}(t +\tau/2 - \mu \tau) + \psi_{\mu} z^{p}(t - \mu \tau) \;\equiv \;\;0 \, ,
\end{equation}
with linear and possibly time dependent operators
\begin{equation}\notag
\phi_{\mu}:\mathcal{P}^{n} \rightarrow \mathcal{I} \quad, \qquad
\psi_{\mu}:\mathcal{P}^{p} \rightarrow \mathcal{I}
\end{equation}
that map into any common linear space ${\;\mathcal{I}}$ wherein $\mathcal{F}$
has its range.
Similar dual equations
${\; \mathcal{F}^{p} [z_{+}^{p}] [z^{n}] \; \equiv \;0 \;}$
determine the linear evolution of the port states during the connection
cycle.\\
Many, if not almost all (viz. all but a finite number) of the
${\phi_{\mu}}$, ${\psi_{\mu}}$ may be zero. Any maximum
${\mu \in \mathbb{N} \cup \{\infty\}}$, such that $\phi_{\mu} \neq 0$
or $\psi_{\mu} \neq 0$, is called the ({\it dynamical\/}) {\it order\/}
of the model equations, and in general equals the order in time of the
integro-differential equations which physically describe the underlying
dynamical problem in terms of smooth fields.

The following statement provides a theoretical means for computing the
reflection map of equation \eqref{4.4} recursively.

\begin{ltheorem}{\label{T:1}}
Let ${\phi_{0} : \mathcal{P}^{n} \rightarrow \mathcal{I}}$
be \emph{bijective} (i.e. one-to-one and onto).
Then $\mathcal{R}$ solves equation \eqref{4.4} on a \emph{finite}
interval $[\, 0, \, T)$,
if and only if for every $t \in [\, 0, \, T)$
\begin{equation}
\begin{split}
\begin{aligned}\notag
\mathcal{R}[ z_{in}^{n} ]( t ) &= z_{out}^{n}( t ) \, = && \\
&= ( - \phi_{0} )^{-1}
\; \sum \nolimits_{\mu = 0}^{\infty} &&\{ \; \; ( \phi_{\mu} \,
+ \, \psi_{\mu} \, nb \, ) \;
z_{in}^{n}( t - \mu \tau ) \, + \\
& &&\! + ( \phi_{\mu + 1} \,
+ \, \psi_{\mu} \, nb \,) \; z_{out}^{n}( t - \tau - \mu \tau )
\; \; \} \;.
\end{aligned}
\end{split}
\end{equation}
\end{ltheorem}

\begin{proof}
Substituting in \eqref{4.4} for $z^{\,p, \,n\;}$ the right-hand sides
of \eqref{3.13} and using equations \eqref{3.15} yields
for ${t < T }\,$ on a time domain trivially extended over the negative
axis
\begin{equation}
\begin{split}
\begin{aligned}\notag
\sum \nolimits_{\mu=0}^{\infty} \; \{ \; \; &\phi_{\mu} \;
[ \, z_{in}^{n}( t + \tau/2 - \mu \tau ) \,
+ \, z_{out}^{n}( t +\tau/2 - \mu \tau ) \, ]
\; + \\ + \; &\psi_{\mu} \, nb \;
[ \, z_{in}^{n}( t + \tau/2 - \mu \tau ) \, 
+ \, z_{out}^{n}( t - \tau/2 - \mu \tau ) \, ]
\; \; \} \, \equiv \, 0 \; ,
\end{aligned}
\end{split}
\end{equation}
which for invertible $\phi_{0}$ and ${ t' :\,= t + \tau/2 }\/$
is equivalent to the recursion formula of the theorem.
\end{proof}

\begin{remark}
For every $t < T < \infty$ the sum in Theorem \ref{T:1} is actually
finite, hence convergence is not a question as long as one abstains
from considering limits. For finite $T\,$, the theorem conveys
\emph{purely algebraic relations} inherited from the 
\emph{structure of the DSC process} as laid down by
(\ref{3.13}, \ref{3.14}) in section~ \ref{S:sec3}.
The question of stability of a process for $T \to \infty$ is being
addressed below.
\end{remark}

If $\phi_{0}$ is \emph{not} bijective, then the model equations
\eqref{4.4} are incomplete in that they do not determine a
well-defined and unique reflection propagator in every cell.
We therefore require that the \emph{uniqueness conditition}
\begin{itemize}
\item[ ] (U)
$\qquad \phi_{0} : \mathcal{P}^{n} \rightarrow \mathcal{I}$
is \emph{bijective} \; \; (\,i.e. a linear isomorphism \,)
\end{itemize}
be always satisfied. Clearly, \eqref{4.4} then defines an explicit
scheme.

The model equations of the classical TLM method discretize Maxwell's
equations. So, they are linear and first order in time, viz.
of the type
\begin{equation}\label{4.5}\centering
\phi_{0} \, z^{n} ( t + \tau / 2 ) \, + \, \phi_{1} \, z^{n} ( t - \tau / 2 ) 
\, + \, \psi_{0} \, z^{p} ( t ) \; \;,
\end{equation}
with time independent operators
$\phi_{\, 0, \, 1}$ and $\psi_{\, 0}$,
which are derived in \cite{He2, He1}, for instance. Of the same type
- ~viz. linear and of first dynamical order~ - are the discretized
diffusion equations.
Theorem \nolinebreak \ref{T:1} applied to this special situation yields

\begin{lcorollary}{\label{C:1}}
The \textnormal{(}for a given time step \emph{unique}\textnormal{)}
reflection map that solves on
a finite interval the first order linear equations \eqref{4.5} with
\emph{time independent} ${\phi_{\, 0, \, 1}}$ and ${\psi_{\, 0}}$ is
\begin{equation}\label{4.6}
\begin{split}
\begin{aligned}
\mathcal{R} [ z_{in}^{n} ] ( t ) \; &= \; z_{out}^{n} ( t ) \; = \\
&= \; K z_{in}^{n} ( t ) \, + \, L \,
\sum \nolimits_{\nu = 1}^{\infty} N^{\nu - 1} M
z_{in}^{n} ( t - \nu \tau ) \; ,
\end{aligned}
\end{split}
\end{equation}
wherein
\begin{equation}
\begin{split}
\begin{aligned}\notag
K \, &= \; - Id \, - \phi_{0}^{-1} \, \psi_{0} \, nb \\
L \; &= \; - \, \phi_{0}^{-1} \\
M \; &= \; + \, \phi_{1} \, + \, ( \, \phi_{1} + \psi_{0} \, nb \, ) \,
\underbrace{( \, - \, \phi_{0}^{-1} \, )
( \, \phi_{0} + \psi_{0} \, nb \, )}_{K} \; \quad \\
N \; &= \; - \, ( \, \phi_{1} + \psi_{0} \, nb \, ) \, \phi_{0}^{-1} \; .
\end{aligned}
\end{split}
\end{equation}
\end{lcorollary}

\begin{proof}
By induction, applying Theorem \ref{T:1} to an incident Dirac pulse
\begin{equation}\notag
z_{in}^{n} ( t ) \, 
:\,= \, \bf{z} \; \chi_{ [ \, 0,\, \tau\, ) \, } ( \, t - \tau/2 \, ) \;
\end{equation}
with any fixed state vector $\bf{z} \in \mathcal{P}^{n}\,$; \,
${{\chi}_{_{I}}}$ denotes the characteristic function of interval $I$
(which equals 1 for every argument in $I$ and 0 elsewhere).
By linearity the statement then holds for arbitrary incident processes
$z_{in}^{n}( t )$, each of them being writable as a superposition of
Dirac processes that start at subsequent time steps.
\end{proof}

\begin{remark}
A sufficient condition for convergence of the propagator series \eqref{4.6}
(applied to any finite incident pulse ${z_{in}^{n} \; }$)
in the limit $t \to \infty$ is obviously
\begin{equation}\label{4.7}\centering
\parallel N \parallel \; \; < \; 1 \quad ,
\end{equation}
where $\,\parallel...\parallel\,$ is any submultiplicative norm,
e.g. the \textsc{Hilbert} (\emph{spectral\,}) norm
\begin{equation}\label{4.8}\centering
\parallel N \parallel_{H} \; \; :\,= \; max \, \{ \, | \, \lambda \, | \,
\text{; $\, \lambda$ eigenvalue of N} \, \}.
\end{equation}
In fact, the latter condition is sufficient for algorithm stability of
the Maxwell field TLM model, as shown in \cite{He1,He2}.
\emph{Any} first order linear process is clearly stable, if
the Hilbert norms of $K, L, M, N$ are bounded by 1 (strictly for $K$
and $N\/$),
since the propagator $\mathcal{R}$ then is \emph{contractive}.
This is in general ensured with bounds for the time step
(if necessary, in combination with a transformation \eqref{4.10}).

For linear model equations of \emph{any\/} finite order, recursion
formulae that generalize \eqref{4.6} are easily derived from Theorem
\nolinebreak \ref{T:1}; e.g. \cite{He5}, equation \nolinebreak (32).
\end{remark}

\begin{lcorollary}{\label{C:2}}
With K, L, M, N as above, the DSC process solving \eqref{4.5}
permits a representation as a 'deflected' scattering process
\begin{equation}\label{4.9}\centering
\left(
\begin{matrix}
z_{out}^{n} ( \, t \, ) \\
d( \, t \, )
\end{matrix}
\right) = 
\left(
\begin{matrix}
K & L \\
M & N
\end{matrix}
\right) \;
\left(
\begin{matrix}
z_{in}^{n} ( \, t \, ) \\
\, d( \, t - \mu \tau \, ) \,
\end{matrix} 
\right) \; ,
\end{equation}
the \emph{deflection\/} ${\,d(t)\,}$ being completely recursively
defined with the \emph{(universally maintained\,)} initial conditions 
${\,d_{|t<0}\,=\,0\,}$.
\end{lcorollary}

\begin{proof}
Immediate consequence of Corollary \ref{C:1}
\end{proof}

\begin{remark}
Note that there remains an arbitrariness in the definition of the
operators $\,L,\,\,M\,,$ and $\,N\,$ in Corollaries \ref{C:1}, \ref{C:2}.
In fact, any invertible transformation
${I \, : \, \mathcal{I} \to \mathcal{I}\,}$ together with the
simultaneous replacements of $\,L,\,M\,,$ and $\,N\,$ by, respectively,
\begin{equation}\label{4.10}\centering
L \sptilde \, :\,= \, L \circ I^{-1} \; , \quad M \sptilde \, :\,= \, I \circ M
\, , \; \text{and}
\quad N \sptilde \, :\,= \, I \circ N \circ I^{-1}
\end{equation}
clearly do not alter the propagator $\mathcal{R}$, i.e. generates
the same DSC process.
\end{remark}

The corollaries offer solutions of first order linear equations with
time independent operator coefficients which are complete in the sense
of (U).
They thus cover the entire field of classical TLM (with connection maps
that are trivial in reducing essentially to identity).

In certain situations it may be useful, or necessary, to integrate some
new (possibly non-linear) interactions into a given DSC model.
Sometimes this can be carried out by adding suitable coupling terms
to the equations of the yet existing model.
We therefore consider \emph{perturbed} model equations of the type
\begin{equation}\label{4.11}
\mathcal{F}^{\,\sptilde}[z_{+}^{n}][z^{p}] \, \equiv \,
\mathcal{F}^{\, n}[z_{+}^{n}][z^{p}] + \mathcal{J}[z_{-}^{n}][z^{p}] \,
\equiv \, 0
\end{equation}
(here exemplary for nodal perturbations), wherein
$\mathcal{J}$ denotes any causal map into $\mathcal{I}$.
The (${-\tau/2}$) time shift in the first argument of $\mathcal{J}$
again synchronizes port and node switching.
Note that the time shift is negative, here.
This is to ensure that the perturbation $\mathcal{J}$ cannot destroy
the uniqueness conditions (U) in the case of a
\emph{linear} function $\mathcal{F}^{\, n}$, and to preserve expliciteness
of the updating relations, in general. The importance of this condition
becomes clear in the proof of the following formula.

\begin{lproposition}[Deflection Formula]\label{P:2}
Let the reflection map $\mathcal{R}$ generate solutions of the linear
equations \eqref{4.4} on a finite interval ${ I = [\,0,\,T)}$.
Then ${\,\mathcal{R}^{\, \sptilde}}$ solves the perturbed equations
\eqref{4.11} on $I$, if and only if the so-called
\emph{deflection}
\begin{equation}\notag
\mathcal{D} \quad :\,= \quad \mathcal{R}^{\, \sptilde} \; - \; \; \mathcal{R}
\end{equation}
\raggedright{satisfies recursively}
\begin{equation}\label{4.13}
\begin{split}
\begin{aligned}[2]\notag
\; \phi_{0\mid t} \, \mathcal{D}_{\mid t + \tau/2} \quad = \quad
&- \, &&\mathcal{J}[z_{-}^{\,n}][z^{\,p}]_{\mid t} \\
&- \, &&\sum \nolimits_{\mu \in \mathbb{N}} \,
( \, \phi_{\mu+1 \mid t} \, + \, \psi_{\mu \mid t} \, nb \, )
\mathcal{D}_{\mid t - \tau/2 - \mu \tau} \quad .
\end{aligned}
\end{split}
\end{equation}
\end{lproposition}

\begin{proof}
By definition, $\,\mathcal{R}^{\, \sptilde}$ solves equations \eqref{4.11}
on $I$, if and only if for every incident sequence $[\, z_{in}^{n}\,]\,$ and
\begin{equation}\notag
\begin{aligned}[1]
z_{out}^{n} \; &:\,= \; \mathcal{R}^{\, \sptilde}[z_{in}^{n}] \quad ,
\quad &&z \; \, :\,= \; z_{in} \, + \, z_{out}
\end{aligned}
\end{equation}
holds
\begin{equation}\notag
\begin{split}
\begin{aligned}[1]
\mathcal{F}^{\, \sptilde}[\, z_{+}^{n}\, ][\, z^{p}\, ] \; &=  \\
&\hspace{-1.2cm}= \, \mathcal{F}^{\, n}[\, z_{+}^{n}\, ][\, z^{p}\, ] \, +
\, \mathcal{J}[\, z_{-}^{n}\,][\, z^{p}] \; = \, 0 \quad .
\end{aligned}
\end{split}
\end{equation}
However, with
\begin{equation}\notag
\begin{aligned}[1]
w_{out}^{n} \; &:\,= \; \; \mathcal{R}[z_{in}^{n}] \quad ,
\quad &&w \; \, :\,= \; z_{in} \, + \, w_{out}
\end{aligned}
\end{equation}
(\,hence \, $\mathcal{D} \; = \; z_{out}^{n} \, - \, w_{out}^{n}$\;) \,
and using \eqref{4.4}, \, this is the case iff
\begin{equation}
\begin{split}
\begin{aligned}\notag
- \, \mathcal{J}[\, z_{-}^{n}\,][\,z^{p}\,]_{\mid t} \; = &\\
\hspace{-1.4cm} = \, \sum \nolimits_{\mu \in \mathbb{N}} \,
&\{ \, \phi_{\mu\mid t} \, ( z_{in\mid t + \tau/2 - \mu \tau }^{n} \,
+ \, \underbrace{w_{out\mid t + \tau/2 - \mu \tau }^{n} \,
+ \, \mathcal{D}_{\mid t + \tau/2 - \mu \tau}}_{
z_{out\mid t + \tau/2 - \mu \tau}^{n}} \, ) \; + \\
+ \; &\psi_{\mu\mid t} \,
nb \, ( z_{in\mid t + \tau/2 - \mu \tau}^{n} \,
+ \, \underbrace{w_{out\mid t - \tau/2 - \mu \tau}^{n} \, 
+ \, \mathcal{D}_{\mid t - \tau/2 - \mu \tau}}_{
z_{out\mid t - \tau/2 - \mu \tau}^{n}} \, ) \, \} \,.
\end{aligned}
\end{split}
\end{equation}
In virtue of the linearity of the $\phi_{\mu}$ , $\psi_{\mu}\,$,
the latter identity holds iff
\begin{equation}
\begin{split}
\begin{aligned}\notag
- \, \mathcal{J}[\, z_{-}^{n}\,][\,z^{p}\,]_{\mid t} \; &= \;
\underbrace{\mathcal{F}^{\, n}[\, w_{+}^{n}\, ][\, w^{p}\, ]}_{
\qquad \quad \, \mathbf{0} \, \text{\, by \eqref{4.4}}} \, + \\
&\quad + \; \phi_{0 \mid t} \, \mathcal{D}_{\mid t + \tau/2} +
\sum \nolimits_{\mu \in \mathbb{N}} \, ( \phi_{\mu+1 \mid t}
+ \psi_{\mu \mid t} \, nb \, )\,
\mathcal{D}_{\mid t - \tau/2 - \mu \tau} \;,
\end{aligned}
\end{split}
\end{equation}
which is the recurrence relations of the proposition.
\end{proof}

\begin{corollary}[Deflected processes]{\label{C:3}}
Let $\mathcal{R}$ solve \eqref{4.4} on a finite interval $I$
and $\phi_{0}:\mathcal{P} \to \mathcal{I}$
be any \emph{bijective} operator
\textnormal{(}that thus satisfies the completeness conditions
\textnormal{(U))}\,.
Then
\begin{equation}
\begin{split}
\begin{aligned}\notag
\mathcal{D}_{\mid t + \tau/2 } \; :\,= \, &- \phi_{0 \mid t}^{-1} 
\; \{ \; \mathcal{J}[\, z_{-}^n \,][\, z^{p}]_{\mid t} \, + \,
\sum \nolimits_{\mu \in \mathbb{N}} (\, \phi_{\mu + 1 \mid t} \, + \,
\psi_{\mu \mid t} \, ) \, \mathcal{D}_{\mid t - \tau/2 - \mu \tau} \; \}
\end{aligned}
\end{split}
\end{equation}
\raggedright{with initial conditions
$\mathcal{D}_{\mid t < 0} \, \equiv \, 0$
defines recursively a causal operator $\mathcal{D}$, such that
$\mathcal{R}^{\, \sptilde} \, :\,= \; \mathcal{R} \, + \, \mathcal{D}$
solves equations \eqref{4.11} on $I$ }.
\end{corollary}

Concluding this section, we stress once again that Theorem \ref{T:1}
and the ensuing propositions and corollaries apply just as well to the
connection cycle, i.e. to cell interface scattering, provided the replacements
\eqref{4.3n} by \eqref{4.3p}, $\mathcal{R}$ by $\mathcal{C}$,
$\mathcal{P}_{\zeta}$ by $\mathcal{P}_{\iota}$, port by node superscripts,
and incoming by outgoing fields (and vice-versa) are simultaneously made.
- \nolinebreak Note, however, that any \emph{excitations}
may \emph{temporarily violate\/} the model equations at a \emph{mesh
boundary\/} face. The model developer is encouraged to care for physically
consistent excitations.

\vspace{-0.3cm}
\section{A heat propagation scheme in non-orthogonal mesh}\label{S:sec5}
The physical interpretation underlying the following application relates
a smoothly varying (viz. in time and space continuously differentiable,
$C^1$-) temperature field $T$,
evaluated as $T^{p}$ at the face centre points and as $T^{n}$
in the nodes of a of non-orthogonal hexahedral mesh, to total states
${ z_{\mu}^{p,n}}$ of a DSC model. 
In fact, we derive the model equations for the connection and reflection
cycles of a DSC heat propagation (diffusion) scheme. 
Since the equations are linear and of dynamical order 0 and 1,
respectively - \nolinebreak as will be seen \nolinebreak - they can be
processed, following the guidelines of the last section.
In the end, we display some computational results of a dispersion test
carried out with this model.

In order to simplify the notation we follow Einstein's convention
to sum up over identical right-hand \nolinebreak (!) sub and
superscripts within all terms where such are present (summation is
not carried out over any index that also appears somewhere at the
left-hand side of a pertinent symbol - \nolinebreak thus, in
${(-1)}^{\kappa} \, a_{\kappa}^{\lambda} \, b_{\lambda} \, \, _{\kappa} c \,$ \,
the sum is made over $\lambda$ but not over $\kappa\,$).
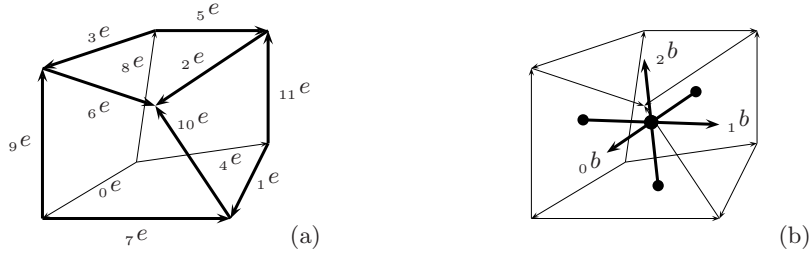
\begin{figure}[!h]\centering
\vspace{-.5cm}
\setlength{\unitlength}{1.cm}
\begin{pspicture}(-1.4,-.5)(20,3.5)\centering
\psset{xunit=.5cm,yunit=.5cm}
\psline[linewidth=0.1mm]{->}(2.5,1.5)(0.0,0.0)
\psline[linewidth=0.1mm]{->}(2.5,1.5)(6.0,2.0)
\psline[linewidth=0.1mm]{->}(2.5,1.5)(3.0,5.0)
\psline[linewidth=0.4mm]{->}(0.0,0.0)(0.0,4.0)
\psline[linewidth=0.4mm]{->}(0.0,0.0)(5.0,0.0)
\psline[linewidth=0.4mm]{->}(6.0,2.0)(6.0,5.0)
\psline[linewidth=0.4mm]{->}(6.0,2.0)(5.0,0.0)
\psline[linewidth=0.4mm]{->}(5.0,0.0)(3.0,3.0)
\psline[linewidth=0.4mm]{->}(3.0,5.0)(0.0,4.0)
\psline[linewidth=0.4mm]{->}(3.0,5.0)(6.0,5.0)
\psline[linewidth=0.4mm]{->}(0.0,4.0)(3.0,3.0)
\psline[linewidth=0.4mm]{->}(6.0,5.0)(3.0,3.0)
\rput(1.8,0.8){$_{_{0}} e$}
\rput(6.0,1.0){$_{_{1}} e$}
\rput(4.0,4.2){$_{_{2}} e$}
\rput(1.5,4.9){$_{_{3}} e$}
\rput(5.0,1.5){$_{_{4}} e$}
\rput(4.4,5.5){$_{_{5}} e$}
\rput(1.5,2.9){$_{_{6}} e$}
\rput(2.5,-.5){$_{_{7}} e$}
\rput(2.4,4.1){$_{_{8}} e$}
\rput(-.6,2.0){$_{_{9}} e$}
\rput(4.0,2.6){$_{_{10}} e$}
\rput(6.7,3.4){$_{_{11}} e$}
\rput(7.0,-0.5){\small{\textnormal{(a)}}}
\psline[linewidth=0.1mm]{->}(15.5,1.5)(13.0,0.0)
\psline[linewidth=0.1mm]{->}(15.5,1.5)(19.0,2.0)
\psline[linewidth=0.1mm]{->}(15.5,1.5)(16.0,5.0)
\psline[linewidth=0.1mm]{->}(13.0,0.0)(13.0,4.0)
\psline[linewidth=0.1mm]{->}(13.0,0.0)(18.0,0.0)
\psline[linewidth=0.1mm]{->}(19.0,2.0)(19.0,5.0)
\psline[linewidth=0.1mm]{->}(19.0,2.0)(18.0,0.0)
\psline[linewidth=0.1mm]{->}(18.0,0.0)(16.0,3.0)
\psline[linewidth=0.1mm]{->}(16.0,5.0)(13.0,4.0)
\psline[linewidth=0.1mm]{->}(16.0,5.0)(19.0,5.0)
\psline[linewidth=0.1mm]{->}(13.0,4.0)(16.0,3.0)
\psline[linewidth=0.1mm]{->}(19.0,5.0)(16.0,3.0)
\psline[showpoints=true,linewidth=0.4mm]{->}(17.375,3.375)(15.0,1.75)
\psline[showpoints=true,linewidth=0.4mm]{->}(14.375,2.625)(18.0,2.50)
\psline[showpoints=true,linewidth=0.4mm]{->}(16.37,.875)(16.0,4.25)
\psline[showpoints=true,
linewidth=0.6mm]{-}(16.1875,2.5625)(16.1875,2.5625) 
\rput(14.55,1.50){$_{_{0}} b$}
\rput(18.5,2.60){$_{_{1}} b$}
\rput(16.6,4.40){$_{_{2}} b$}
\rput(20,-0.5){\small{\textnormal{(b)}}}
\end{pspicture}
\caption{\textsl{Non-orthogonal hexahedral mesh cell. \newline
\textnormal{(a)} Edge vectors. \qquad \qquad 
\textnormal{(b)} Node vectors.}}\label{F:2}
\end{figure}

\begin{figure}[!h]\centering
\vspace{-1.5cm}
\setlength{\unitlength}{1.cm}
\begin{pspicture}(-2.0,-2)(10,4.0)\centering
\psset{xunit=.5cm,yunit=.5cm}
\psline[linewidth=0.2mm]{->}(0.0,4.0)(0.0,1.0)
\psline[linewidth=0.2mm]{->}(0.0,1.0)(5.0,0.0)
\psline[linewidth=0.2mm]{->}(0.0,4.0)(4.0,5.0)
\psline[linewidth=0.2mm]{->}(4.0,5.0)(5.0,0.0)
\psline[showpoints=true,linewidth=0.6mm]{-}(2.25,2.5)(2.25,2.5) 
\psline[linewidth=0.4mm]{->}(4.5,2.5)(0,2.5)
\psline[linewidth=0.4mm]{->}(2.0,4.5)(2.5,0.5)
\psline[showpoints=true,
linewidth=0.4mm]{->}(4.5,2.5)(8.751,3.400) 
\psline[showpoints=true,linewidth=0.4mm]{->}(0.0,2.5)(-1.5,2.5)
\psline[showpoints=true,linewidth=0.4mm]{->}(2.0,4.5)(1.287,7.352)
\psline[showpoints=true,linewidth=0.4mm]{->}(2.5,0.5)(1.650,-3.751)
\rput(1.0,3.0){$_{_{0}}b$}
\rput(3.0,1.5){$_{_{1}}b$}
\rput(6.6,3.6){$_{_{0}}f$}
\rput(-1.2,3.2){$_{_{1}}f$}
\rput(2.2,6.4){$_{_{2}}f$}
\rput(2.7,-1.4){$_{_{3}}f$}
\rput(5.4,2.0){$_{_{0}}T^{p}$}
\rput(-0.8,1.9){$_{_{1}}T^{p}$}
\rput(1.0,4.7){$_{_{2}}T^{p}$}
\rput(1.6,0.0){$_{_{3}}T^{p}$}
\rput(2.9,3.0){$T^{n}$}
\end{pspicture}
\caption{\textsl{Face vectors and
temperature points (nodal section).\qquad\qquad\qquad\qquad\;\;}}\label{F:3}
\vspace{-0.2cm}
\end{figure}
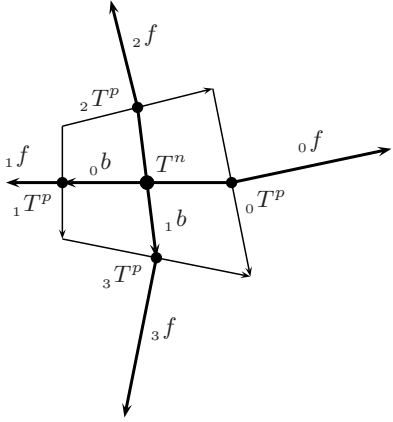

The shape of a hexahedral cell is completely determined by its 12
\emph{edge vectors} ${(_{\nu} e)_{\nu=0,...,11}}\/$.
Also, with the labelling scheme of fig \nolinebreak \ref{F:2}a,
\emph{node
vectors\/} ${(_{\mu} b)_{\mu = 0,1,2}}$ and \emph{face vectors\/}
${(_{\iota} f)_{\iota = 0,...,5}}$ are defined as
\begin{equation}\label{5.1}\centering
\begin{split}
\begin{aligned}
_{\mu} b \quad &:\,= \quad \quad \frac{1}{4}
&& \! \! \sum \nolimits_{\nu = 0}^{3} \, _{_{(4 \mu + \nu)}} e \, \,
&&\mu = 0,1,2 \\
\text{and} \qquad
_{\iota} f \quad &:\,= \quad \, \frac{(-1)^{\iota}}{4}
&& \, ( \, \, _{_{(8 + 2\iota)}} e \,
+ _{_{(9 + 2(\iota + (-1)^{\iota}))}} e \, )
\, \, \land && \\
& &&\; \; \land \, ( \, _{_{(4 + 2\iota)}} e \, 
+ _{_{(5 + 2\iota)}} e \, ) \, \, &&\iota = 0,...,5  \,
\end{aligned}
\end{split}
\end{equation}
(edge vector indices cyclic modulo 12\, and
$\land$ denoting the wedge ('cross') product in $\mathbb{R}^3$).

At every face ${\iota \in \{0,...,5\}}\/$ of a mesh cell and for any 
given $\tau \in \mathbb{R}_{+}\,$ the following time shifted finite
temperature differences in directions ${ _{\mu} b }$ (\,$\mu = 0,1,2\,$)
form a vector valued function
\begin{equation}\label{5.2}\centering
\begin{split}
_{\iota} \! {\nabla}^{B} T_{\mu}\,( \,t\,,\, \tau\, ) \quad : \, = \quad
\begin{cases}
\, 2 \, (-1)^{\iota} (\, T^{n} \, _{\mid \, t - \tau/2} 
- \, _{\iota} T^{p} \, _{\mid \, t} \,) \qquad
&\text{if $\mu \, = \, [ \iota / 2 ]$} \\
\, (\,\, _{2\mu + 1} T^{p} \, - \, _{2\mu} T^{p} \,\,)
\, _{\mid \, t - \tau \, } \qquad
&\text{if $\mu \, \neq \, [ \iota / 2 ]$} \,
\end{cases}
\end{split}
\end{equation}
($\,[\, x \,]$ denotes the integer part of $x \in \mathbb{R}\,$).
The time increments are chosen to attain technical consistence
with the updating conventions of DSC schemes. They do not destroy
convergence, as easily seen:
In fact, in the centre point of face ${\, \iota \,}$ the vector
${\,_{\iota} \! {\nabla}^{B} T\,}$
approximates in the first order of the time increment ${\, \tau \,}$,
and of the linear cell extension, the scalar products of the node
vectors with the temperature gradient ${\, \nabla T \,}$.
Let, precisely, for a fixed centre point on face ${\, \iota \,}$
and $\,\epsilon \in \mathbb{R}_{+}\,$ the \emph{$\epsilon$-scaled cell}
have edge vectors
$\, _{\iota} e\sptilde \, : \, = \, \epsilon \, \, _{\iota} e \,$. 
Let also $\, _{\iota} {\nabla}^{B\sptilde} T_{\mu} \,$ denote
function \eqref{5.2} for the $\epsilon$-scaled cell (with node vectors 
$\, _{\mu} b\sptilde \, = \, \epsilon \, _{\mu} b \;$).
Then at the fixed point holds
\begin{equation}\label{5.3}\centering
\begin{split}
< \, _{\mu} b \, , \, \text{grad($T$)} \, > \, \, \,
= \, \, _{\mu} b \cdot \nabla T \, \,
= \, \, \lim_{\epsilon \to 0} \, \, \lim_{\tau \to 0} \, \,
\frac{1}{\epsilon} \, _{\iota} \! {\nabla}^{B^{\sptilde}} T_{\mu} \, ,
\end{split}
\end{equation}
as immediately follows from the required $C^1$-smoothness of the
temperature field $T$.

To recover, in the same sense and order of approximation, the gradient
${{\nabla} T \,}$ from \eqref{5.2}, observe that for every orthonormal
basis $( _{\nu} u )_{\nu = 0,...,m-1}\,$ of
$V = \mathbb{R}^{m} \, \text{or} \, \, \mathbb{C}^{m}\,$, and for an
arbitrary basis $( _{\mu} b )_{\mu = 0,...,m-1}$ with coordinate
matrix ${\beta_{\nu}^{\mu}} \, = \, {< \, _{\nu} u \, , \, _{\mu} b \, >}$,
the scalar products of any vector $a \in V$ with ${_{\mu} b}$ are
\begin{equation}\label{5.4}\centering
\underbrace{< \, _{\mu} b \, , \, a \, >}_{ \qquad = \, : 
\, \, {\alpha}_{\mu}^{B}} \, \, = \, \sum \nolimits_{\nu = 0}^{m-1} \,
\underbrace{< \, _{\mu} b \, , \, _{\nu} u \, >}_{ \, \, \,
( {\bar{\beta}}_{\mu}^{\nu} ) \, = \, ( {\beta}_{\nu}^{\mu} )^{^{*}}} \,
\underbrace{< \, _{\nu} u \, , \, a \, >}_{ \qquad = \, : \, \, {\alpha}_{\nu}}
\, \, = \, \bar{\beta}_{\mu}^{\nu} \, {\alpha}_{\nu} \; .
\end{equation}
(At the right-hand side, and henceforth,
we follow Einstein's convention to sum up over identical sub and superscripts
within terms where such are present). It follows that
\begin{equation}\label{5.5}\centering
{\alpha}_{\nu} \, = \,
{\gamma}_{\nu}^{\mu} \alpha_{\mu}^{B} \, ,
\qquad \text{with} \qquad ( {\gamma}_{\nu}^{\mu} ) 
\, : \, = \, {({(\beta_{\nu}^{\mu})}^{*})}^{-1} \quad .
\end{equation}

Loosely speaking, the scalar products of any vector with the basis vectors
${_{\mu} b\,}$ transform into the coordinates of that vector with
respect to an orthonormal basis ${_{\nu} u \,}$ by multiplication
with matrix ${\gamma = (\beta^*)^{-1}\,}$, where
${\beta_{\nu}^{\mu}} \, \, = \, {< \, _{\nu} u \, , \, _{\mu} b \, >}\,$,
i.e. $\beta$ is the matrix of the coordinate (column) vectors
${_{\mu} b \,}$
with respect to the given ON-basis ${_{\nu} u \,}$, and $\gamma$ is the
adjoint inverse of $\beta$.

This applied to the node vector basis ${ _{\mu} b \,}$ and \eqref{5.3}
yields the approximate temperature gradient at face $\iota$ as
\begin{equation}\label{5.6}\centering
_{\iota} \! \nabla T_{\nu} \quad
= \quad {\gamma}_{\nu}^{\mu} \, \, \, _{\iota} \! {\nabla}^{B} T_{\mu}.
\end{equation}

Thus, the heat current \emph{into} the cell through face $\iota$ with
face vector components
${ _{\iota} f^{\nu}} \,
= \, \, {< \, _{\iota} f \, , \, _{\nu} u \, >}\,$, \,
$\nu \in \{ 0,1,2 \}\,$, is
\begin{equation}\label{5.7}\centering
\begin{split}
\begin{aligned}
\qquad _{\iota} J \, \, 
&= \, \, {\lambda}_{H} \, _{\iota} f \, \cdot \, _{\iota} \! \nabla T \, \,
= \, \, \underbrace{{\lambda}_{H} \, _{\iota} f^{\nu} \, \,
{\gamma}_{\nu}^{\mu}}_{ \qquad = \,: \, \, _{\iota} s^{\mu}} \, \,
_{\iota} \! {\nabla}^{B} T_{\mu} \, \,
= \, \  _{\iota} s^{\mu} \, \, _{\iota} \! {\nabla}^{B} T_{\mu} \quad ,
\end{aligned}
\end{split}
\end{equation}
${{\lambda}_{H}}\,$ denoting the heat conductivity in the cell.

The heat current through every interface is conserved, i.e. between any
two adjacent cells $\zeta$, $\chi$ with the common face labelled $\iota$
in cell $\zeta$ and $\kappa$ in $\chi$ applies
\begin{equation}\label{5.8}\centering
_{\iota}^{^{\zeta}} \! J \quad = \quad - \,\, _{\kappa}^{^{\chi}} \! J \, .
\end{equation}
Also, the nodal temperature change in cell $\zeta$ is
\begin{equation}\label{5.9}\centering
\frac{ d \, {T}^{n}}{dt} \quad 
= \quad \frac{1}{c _{v} \, V } \,
( \, \, S \, + \, \sum \nolimits_{\iota = 0}^{5} \,
_{\iota}^{^{\zeta}} \! J \, \, ) \, ,
\end{equation}
where ${c_{v}}$ denotes the heat capacity (per volume), $V$ the cell
volume, and $S$ any heat source(s) in the cell. \\
We finally introduce quantities
${_{\iota} z _{\mu}^{p,n}}$ \, ($\iota \, = \, 0,...,5;$ \, $\mu \, = \,
0,1,2 \,$), which still smoothly vary in time with the
temperature $T$ (and that are hence not yet DSC states, but will
later be updated as such)
\begin{equation}\label{5.10}\centering
\begin{split}
_{\iota} z _{\mu}^{n} \, ( \, t \, ) \quad : \, = \quad
\begin{cases}
\, \, 2 \, (-1)^{\iota} \, \, T^{n} \, _{\mid \, t} \qquad
&\text{if $\mu \, = \, [ \iota / 2 ]$} \, \, \\
\, \, ( \, _{2\mu + 1} T^{p} 
- \, _{2\mu} T^{p} \,)_{\mid \, t - \tau/2} \qquad
&\text{else} \,
\end{cases}
\end{split}
\end{equation}
and
\begin{equation}\label{5.11}
\begin{split}
_{\iota} z _{\mu}^{p} \, ( \, t \, ) \quad : \, = \quad
\begin{cases}
\, \, 2 \, (-1)^{\iota} \, \, _{\iota} T^{p} \, _{\mid \, t} \qquad
&\qquad \qquad \text{if $\mu \, = \, [ \iota / 2 ]$} \\
\, \, \, _{\iota} z _{\mu}^{n} \, ( \, t \, - \, \tau/2 \, ) \qquad
&\qquad \qquad \text{else} \, .
\end{cases}
\end{split}
\end{equation}
From (\ref{5.2}, \ref{5.7}) follows
\begin{equation}\label{5.12}
\begin{split}
_{\iota} J\, _{\mid \, t + \tau / 2} \quad
&= \quad \, _{\iota} s^{\mu}\, ( \, _{\iota} z_{\mu}^{n}\, _{\mid \, t} \,
- \, 2 \, {(-1)}^{\iota} {\delta}_{\mu}^{[\iota/2]} \, \,
_{\iota} T^{p} \,_{\mid \, t + \tau/2}\, ) \\
&= \quad _{\iota} s^{\mu}\, ( \, _{\iota} z_{\mu}^{n}\, _{\mid \, t} \,
- \, {\delta}_{\mu}^{[\iota/2]} \, \, 
_{\iota} z_{\mu}^{p}\, _{\mid \, t + \tau/2} \, )
\quad .
\end{split}
\end{equation}
This, with \eqref{5.8} and continuity of the temperature at the interface,
\, $_{\iota}^{^{\zeta}} T \, ^{p} \, = \, _{\kappa}^{^{\chi}} T \, ^{p}\,$, \,
together with \eqref{5.11} imply
\begin{equation}\label{5.13}
\begin{split}
_{\iota}^{^{\zeta}} z \, _{\mu}^{p} \, _{\mid \, t + \tau/2} \, \, = \,
\begin{cases}
\, \frac{\, _{\iota}^{^{\zeta}} s \, ^{\mu} \, \, 
_{\iota}^{^{\zeta}} z \, _{\mu}^{n} \, _{\mid \, t} \,
+ \, _{\kappa}^{^{\chi}} s \, ^{\nu} \, \, \,
_{\kappa}^{^{\chi}} z \, _{\nu}^{n} \, _{\mid \, t}} 
{_{\iota}^{^{\zeta}} s \, ^{[\iota/2]}
+ \, ( -1 )^{\iota + \kappa} \, _{\kappa}^{^{\chi}} s \, ^{[\kappa /2 ]}}
\qquad \, \,
&\text{if $\mu = [\iota/2]$} \vspace{0.2cm} \\
\, \, _{\iota}^{^{\zeta}} z \, _{\mu}^{n} \, _{\mid \, t}
\qquad \, \,
&\text{else} \, ,
\end{cases}
\end{split}
\end{equation}
which form a complete set of recurrence relations for ${z^{p}}$
(\,given ${z^{n}}\/$) and so can be taken as model equations for the
connection cycle of a DSC algorithm.

Equations \eqref{5.9} discretely integrated in the time balanced form with
increment $\tau$ yield
\begin{equation}\label{5.14}
\begin{split}
T^{n} \, _{\mid \, t \, + \tau/2 } \quad
&= \quad
T^{n} \, _{\mid \, t \, - \tau/2 } \,
+ \frac{\tau}{c_{v} \, V} \, ( \, \, S \, + \, \, \sum \nolimits_{\iota = 0}^{5}
\, _{\iota} J \, _{\mid t } \, \, ),
\end{split}
\end{equation}
i.e., with (\ref{5.10}, \ref{5.11}, \ref{5.12}), the recurrence relations
\begin{equation}
\begin{split}\label{5.15}
_{\iota} z\, _{\mu}^{n} \, _{\mid \, t + \tau/2 }\,=\,
\begin{cases}
\,_{\iota} z\,_{[\iota / 2]}^{n}\,_{\mid\,t - \tau / 2}\,
+\frac{{(-1)}^{\iota}\,\tau}{2\,c_{v}\,V}\,\,\{\,\,S\,\,
+\\
\quad\;+\,\,\sum\nolimits_{\kappa = 0}^{5}\,
_{\kappa} s^{\nu}\,(\, _{\kappa} z\,_{\nu}^{n}\,_{\mid\,t - \tau / 2}\,
-\,{\delta}_{\nu}^{[\kappa/2]}\,
_{\kappa} z \,_{\nu}^{p}\,_{\mid\,t}\,)\,\}
&\text{if $\mu = [ \iota / 2 ]$} \vspace{.2cm}\\
\,-\,\frac{1}{2}\,(\,_{2 \mu + 1 } z\,_{\mu}^{p}\,+
_{2\mu} z \,_{\mu}^{p}\,)\,_{\mid\,t}
&\text{else}\,,
\end{cases}
\end{split}
\end{equation}
which provide a complete set of model equations for the reflection cycle
of a DSC algorithm. Note that the first line, modulo the factor
$2\,(-1)^{\iota}\,$, always updates the nodal temperature. Of course,
this has to be carried out only once per cell and iteration cycle.
In this - ~typical ~- example the dual state space concept of DSC
(\,needed by Johns' cycle\,) creates a redundancy, which is a
price for process parallizability within either parts of the 
connection-reflection cycle.

Equations (\ref{5.13}, \ref{5.15}) can be directly taken as updating
relations for total quantities of a DSC scheme. Alternatively, they 
may be further processed in deriving reflection and connection maps,
along with stability bounds for the time step.
The proceeding is canonical and amounts in essence to a straightforward
transcription of the model equations along the lines of 
Theorem~\ref{T:1} and corollaries.

It is particularly easy to couple this heat conduction model - ~within
one and the same mesh ~- to a Maxwell field TLM model in the non-orthogonal
setting \cite{He2}. In fact, with the node vector definition in \cite{He4},
the total node voltages are just the scalar products of $_{\mu} b$ with the
electric field, hence the dielectric losses and heat sources per cell are
\begin{equation}\centering\label{5.16}
S \quad = \quad \frac{1}{2} \,
\sigma \, V \, E^{\nu} \, \overline{E_{\nu}} \quad
= \quad \frac{1}{2} \, \sigma \, V \, \sum \nolimits_{\nu} \, \left|
{\gamma}_{\nu}^{\mu} \, U_{\mu}^{n} \right|^{2} \; ,
\end{equation}
for a frequency domain (complex) TLM algorithm, cf. \cite{He5}\, ;
$\sigma = 2 \pi f \, \epsilon \, \text{tan}( \delta )$
denotes the effective loss current conductivity at frequency $f$
in a mesh cell of absolute permittivity $\epsilon$ and dielectric loss
factor $\text{tan}( \delta )\,$; $\gamma \, = \, ( \beta^{*} )^{-1}\,$
as in \nolinebreak \eqref{5.5}. In \textsc{Spinner}'s Maxwell field
solver the model couples, in addition, to magnetic and skin effect losses.

\begin{figure}[!h]\centering
\vspace{-1.5cm}
\setlength{\unitlength}{1.cm}
\begin{pspicture}(0.0,0.0)(13.0,4.5)\centering
\psset{xunit=1.0cm,yunit=1.0cm}
\rput(1.7,1.67){\includegraphics[scale=0.3150,clip=0]{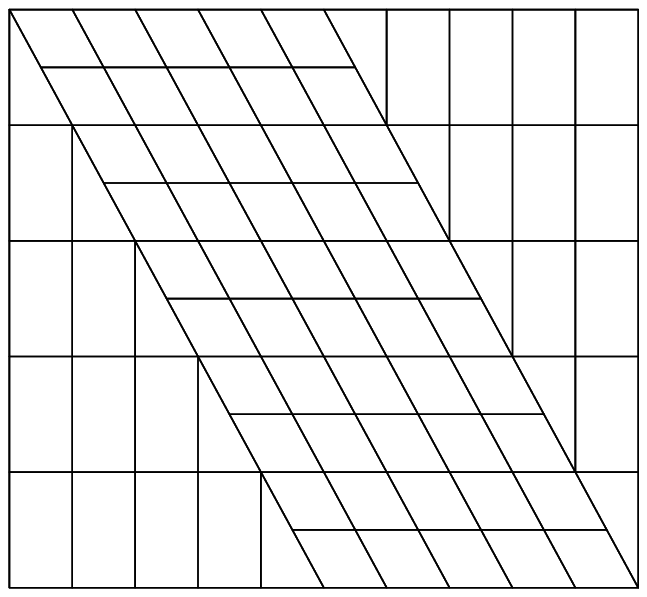}}
\rput(0.25,0.75){\small (a)}
\rput(5.8,1.45){\includegraphics[scale=0.3000,clip=0]{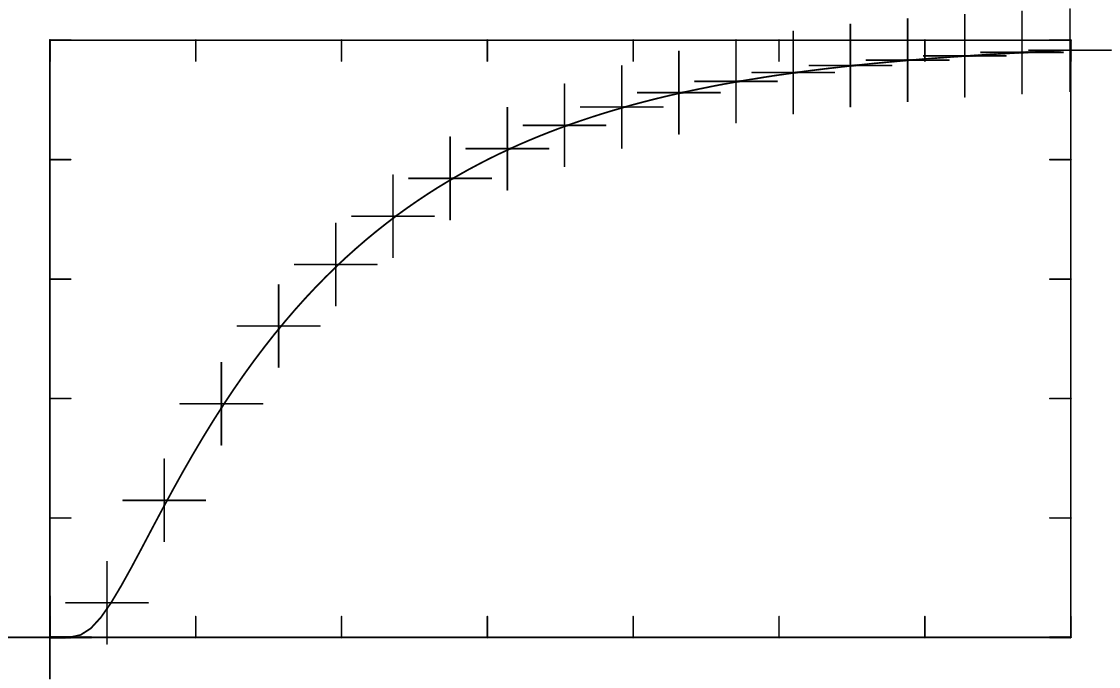}}
\rput(4.0,0.75){\small (b)}
\rput(10.6,1.45){\includegraphics[scale=0.3000,clip=0]{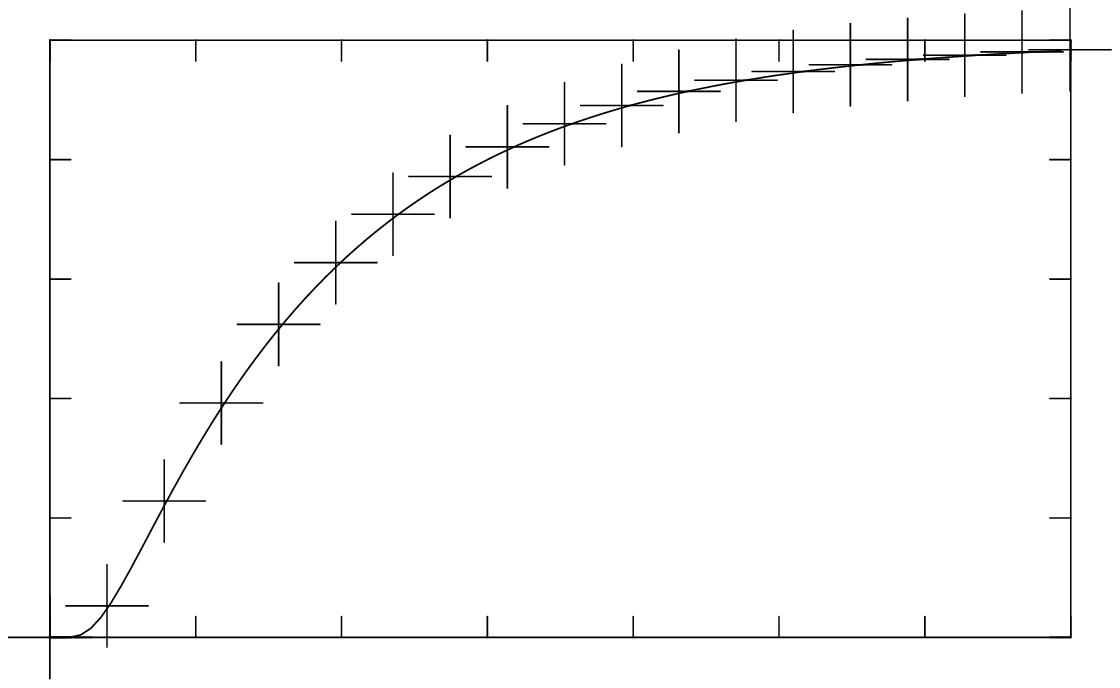}}
\rput(8.8,0.75){\small (c)}
\end{pspicture}
\caption{\textsl{Transverse heat conduction over a square mesh using
non-orthogonal cells. A Heaviside temperature step is imposed on one
side and the transient temperature computed at the opposite side,
assuming adiabatic boundary conditions on all but the heated sides.
\newline
DSC results \textnormal{(+)} are plotted over analytical solution
(curve).
\newline
\textnormal{(a)} The mesh. \, 
\textnormal{(b)} Horizontal \quad
\textnormal{(c)} vertical propagation.}}\label{F:4}
\end{figure}

Finally, Fig \nolinebreak \ref{F:4} displays the result of a dispersion
test, computed in a square mesh using non-orthogonal cells.
It turns out that the heat conduction properties of the mesh are highly
insensitive to cell shape and orientation (as of course should be the case).

\vspace{-0.2cm}
\section{Conclusions}\label{C:sec6}
Johns' TLM algorithm can be extended with benefit in two major directions
by replacing transmission line links between cells with abstract scattering
channels in terms of paired distributions and in admitting non-trivial cell
interface scattering.
Executing this program lead us in this paper to a new class of Dual Scattering
Channel schemes which offer enhanced modeling potentiality and canonical
techniques for stable algorithm design.\\
\textsc{SPINNER}'s implementation of a heat propagation scheme coupled to
a lossy Maxwell field illustrates the approach.\\
The connection and reflection cycles of a DSC process are (either) completely
parallelizable, which can be turned into account in computational performance.
DSC schemes open a challenging field to future research. Applications to fluid
dynamics are presently under examination.
\vspace{-0.2cm}

\textsc{Spinner} RF Lab, Erzgiessereistr. 30, DE-80335 M\"unchen, Germany \\
E-mail address:\; s.hein@spinner.de 
\end{document}